\documentclass[12pt]{article}

\usepackage{amssymb}
\usepackage{amsmath}
\usepackage{amsthm}
\usepackage{amsfonts}

\newcommand{\rmK}{{\rm K}}

\newcommand{\mmp}{\mathbb{P}}

\newcommand{\me}{\mathbb{E}}
\newcommand{\mr}{\mathbb{R}}
\newcommand{\mn}{\mathbb{N}}

\newcommand{\lin}{\underset{n\to\infty}{\lim}}

\newtheorem{thm}{Theorem}[section]
\newtheorem{lemma}[thm]{Lemma}

\newtheorem{assertion}[thm]{Proposition}
\newtheorem{example}[thm]{Example}
\theoremstyle{definition}

\theoremstyle{remark}

\begin{document}
\title{On asymptotics of exchangeable coalescents with multiple collisions}
\date{\today}

\author{Alexander Gnedin\footnote{
Department of Mathematics, Utrecht University, Postbus 80010, 3508
TA Utrecht, The Netherlands, e-mail: A.V.Gnedin@uu.nl},\quad Alex
Iksanov\footnote{ Faculty of Cybernetics, National T. Shevchenko
University of Kiev, 01033 Kiev, Ukraine,\newline e-mail:
iksan@unicyb.kiev.ua}\quad and\quad Martin M\"ohle\footnote{
Mathematical Institute, University of D\"usseldorf,
Universit\"atsstra{\ss}e 1, 40225 D\"usseldorf, Germany, e-mail:
moehle@math.uni-duesseldorf.de}}

\maketitle
\begin{abstract}
\noindent We study the number of collisions $X_n$ of an
exchangeable coalescent with multiple collisions
($\Lambda$-coalescent) which starts with $n$ particles and is
driven by rates determined by a finite characteristic measure
$\nu({\rm d}x)=x^{-2}\Lambda({\rm d}x)$. Via a coupling technique
we derive limiting laws of $X_n$, using previous results on
regenerative compositions derived from stick-breaking partitions
of the unit interval. The possible limiting laws of $X_n$ include
normal, stable with index $1\le\alpha<2$ and Mittag-Leffler
distributions. The results apply, in particular, to the case when
$\nu$ is a beta$(a-2,b)$ distribution with parameters $a>2$ and
$b>0$. The approach taken allows to derive asymptotics of three
other functionals of the coalescent, the absorption time, the
length of an external branch chosen at random from the $n$
external branches, and the number of collision events that occur
before the randomly selected external branch coalesces with one of
its neighbours.
\end{abstract}

\section{Introduction} \label{intro}
Let $\Pi_n=(\Pi_n(t),t\geq 0)$ be an exchangeable  coalescent with
multiple collisions, as introduced by Pitman \cite{Pit} and
Sagitov \cite{Sag}. This is a c{\'a}dl{\'a}g Markov process that
starts at $t=0$ with $n$ particles and evolves according to the
following dynamics. When $m$ particles are present, each $k$-tuple
of them collides and merges to form a single particle at rate
\begin{equation}\label{lrates}
\lambda_{m,k}=\int_{[0,1]} x^k(1-x)^{m-k}\nu({\rm d}x), \qquad
2\leq k\leq m
\end{equation}
The definition of the rates
requires that the {\it characteristic measure} $\nu$ had a finite
second moment. These processes are sometimes called
`lambda-coalescents', following Pitman's  paper \cite{Pit} where
the rates were introduced by means of the finite measure
$\Lambda({\rm d}x)=x^2\nu({\rm d}x)$.

\par
The coalescent starting with $n$ particles is usually understood
as a process with values in the set of partitions of
$[n]=\{1,\ldots,n\}$. To that end, one needs to label the initial
particles by integers and, tracking their evolution, to record
which of the initial particles comprise a merger. In this context,
the exchangeability means that the restriction of the process to
$[m]$ (or any other subset of $[n]$ with $m<n$ elements) is a
stochastic copy of the process starting with $m$ particles,
 see \cite{Bertoin, CSP} for survey and details.
For most of our purposes, however, it is sufficient to identify
 the coalescent with a pure death process of counting the particles,
with state-space $[n]$ and transition rates
\begin{equation*}\label{rate-g}
g_{n,m}={n\choose m-1}\lambda_{n,n-m+1},~~~~~~1\leq m\leq n-1
\end{equation*}
so that the total transition rate from state $n$ is
\begin{equation*}\label{totalrate}
g_n=\sum_{m=1}^{n-1}g_{n,m}= \int_{[0,1]}
(1-(1-x)^n-nx(1-x)^{n-1})\nu({\rm d}x).
\end{equation*}
In terms of the richer partition-valued process, the event
$\Pi_n(t)=m$ occurs when  at time $t$ the partition of $[n]$ has
$m$ blocks, and a merge of some $k\geq 2$ blocks into a single
block corresponds to a transition from $m$ to $m-k+1$.

\par The coalescent starts in the state $\Pi_n(0)=n$ and terminates in
the state $1$, as all particles eventually merge in one. The {\it
number of collisions} $X_n$ is the number of transitions which the
process needs to proceed from $n$ to the terminal state $1$.
 The {\it absorption time}
$\tau_n=\inf\{t\ge 0:\Pi_n(t)=1\}$ is the time when the process
enters the terminal state. Note that in biological applications of
the coalescent theory the absorption time $\tau_n$ has the meaning
of time back to the most recent common ancestor. The main focus of
this paper is on the weak asymptotic behaviour of $X_n$ and
$\tau_n$ as $n$ tends to infinity. In particular, the main result
(Theorem \ref{basic1}) is a criterion describing all possible
limit laws of, properly normalized and centered, $X_n$: these can
be normal, stable with index $1\le\alpha<2$ and Mittag-Leffler
distributions. Another result (Theorem \ref{main1}) provides
sufficient conditions under which $\tau_n$ possesses the same
collection of limit laws. The limiting results for both, $X_n$ and
$\tau_n$, will be obtained by comparison with a simpler process
which was partly investigated in connection with regenerative
composition structures \cite{Gne, GINR}.

\par Intuitively, a high concentration of the measure $\nu$ near
$0$ increases the intensity of small jumps of $\Pi_n$, hence
increases the number of collisions $X_n$. This is well exemplified
by different kinds of the asymptotic behaviour of $X_n$  for
measures $\nu$ with densities $\nu({\rm
d}x)=c\,x^{a-3}(1-x)^{b-1}{\rm d}x$, with $a,b,c >0$. For $0<a<1$
and $b>0$ the order of growth of $X_n$ is linear \cite{Delmas,
GneYak, IksMoe2}, for $a=1$ and $b=1$ the order is $n/\log n$
\cite{DrmIksMoeRoe, DrmIksMoeRoe2, IksMoe2, Pan}, for $1<a<2$ and
$b=1$ the order is $n^{2-a}$ \cite{IksMoe2}, while for $a=2$ and
$b>0$ the order is $\log^2n$ \cite{IksMarMoe}. The class of beta
distributions covered by the present paper corresponds to the
range $a>2$, when $X_n$ exhibits a logarithmic growth.
See \cite{BBS} for other aspects of the beta-coalescents.

\par We shall assume throughout that the characteristic measure
$\nu$ is a probability measure on the open interval $(0,1)$ such
that
\begin{itemize}
\item[(i)] the support of $\nu$ is not contained in a sequence of
the form $(1-\delta\gamma^n, n=0,1,\ldots)$ for some $\delta>0$
and $\gamma\in (0,1)$, \item[(ii)] $\nu$ satisfies
\begin{equation}\label{nuat1}
\int_{(0,1)} |\log x|\nu({\rm d}x)<\infty.
\end{equation}
\end{itemize}
Note that the more general case of a finite measure $\nu$ can be
reduced to the case of a probability measure by a linear time
change of the coalescent.

\section{Annihilator and a coupling}

Adopting Pitman \cite[p. 1884]{Pit}, the coalescent process
$\Pi_n$ can be constructed as follows. Let $0<s_1<s_2<\cdots$ be
the arrival times of a unit Poisson process  and let
$1-\eta_1,1-\eta_2,\ldots$ be an iid sample from $\nu$, also
independent of the Poisson process. At time $s_1$ each of the $n$
particles is marked `head' with probability $1-\eta_1$ or `tail'
with probability $\eta_1$. Those particles marked `head' are
merged in one particle, provided there are at least two `heads'.
Thus $k\geq 2$ `heads' imply that $n-k+1$ particles remain after
merging, i.e. $\Pi_n(s_1)=n-k+1$, and $k=1$ or $=0$ `heads' imply
that $\Pi_n(s_1)=n$. The procedure is then iterated at time $s_2$
with the remaining particles and the new probability for `heads'
being $\eta_2$, and so on. It is immediately checked that this
description agrees with the rates (\ref{lrates}). The rules define
the $\Pi_n$'s consistently for all values of $n$, a property used
by Pitman \cite{Pit} to introduce a partition-valued coalescent
process with infinitely many particles.

\par
There is a simpler process on $\{0\}\cup [n]$, whose transition
rate from $m$ to  $m-k$ for $k=1,\ldots,m$ is ${n\choose
m}\lambda_{n,k}$ (where $\lambda_{n,1}$ is defined by the same
integral   (\ref{lrates}) with $k=1$). Let $\rmK_n=(\rmK_n(t),
t\geq 0)$ denote the process starting at $n$ and terminating at
$0$. We shall call this process {\it annihilator} in view of the
following interpretation: when there are $m$ particles, any
$k$-tuple of them collides and annihilates at rate
$\lambda_{m,k}$, so that after the collision $m-k$ particles
remain (we stress that the case $k=1$ of 'singular' collision is
not excluded). A sequential construction of the annihilator is
almost the same as that of the coalescent. At the generic time
$s_j$ of a unit Poisson process  a random value $1-\eta_j$ is
sampled from $\nu$, and each of the remaining particles (if any)
is marked `head' with probability $1-\eta_j$ or `tail' with
probability $\eta_j$, then the particles marked `heads' are
removed.

\par The sequence of positive decrements of $\rmK_n$ is  a random
{\it composition} (ordered partition) of the integer $n$. For
instance, the path $10,7,2,1,0,\ldots$  of $\rmK_{10}$ corresponds
to the composition $(3,5,1,1)$ of the integer $10$. As $n$ varies,
the random compositions are consistent with respect to a sampling
operation, i.e. they comprise a {\it composition structure}, see
\cite{RCS, GDP}. Two features differ the transition rules of
$\Pi_n$ and $\rmK_n$:
\begin{itemize}
\item[(i)] if at time $s_j$ exactly $k\geq 2$ of the remaining particles are marked
`heads' then $\Pi_n$ decrements by $k-1$, while $\rmK_n$ decrements by $k$,
\item[(ii)] if $k=1$ then $\Pi_n$ has no transition at $s_j$, while $\rmK_n$ decrements by $1$.
\end{itemize}
Note that neither $\Pi_n$ nor $\rmK_n$ jump at time $s_j$ if none
of the remaining particles is marked `head'.

\par To couple the coalescent and the annihilator we just
define them using the same unit Poisson process with arrival times
$s_1,s_2,\ldots$ and the same sample $1-\eta_1, 1-\eta_2,\ldots$
from $\nu$. We shall call the initial $n$ particles {\it primary}
and their followers resulting from mergers {\it secondary}. The
primary particles are labelled $1,2,\ldots,n$. If at some time
exactly $k\geq 1$ particles are marked `head', these $k$ particles
are replaced by one secondary particle. The variable $\rmK_n(t)$
is identified with the number of primary particles alive at time
$t$, and  $\Pi_n(t)$ is identified with the total number of
primary and secondary particles at time $t$ (thus by this coupling
we always have $\Pi_n(t)\geq \rmK_n(t)$). For instance, if at time
$s_1$ the particles $2,3$ and $6$ are marked `head', the number of
primary particles decreases by $3$, and one secondary particle is
born.

\par Let $K_n$ be the number of transitions of $\rmK_n$ as the
process proceeds from $n$ to $0$, and let $K_{n,1}$ be the number
of jumps of size $1$ (singular collisions). Introduce the
absorption time of the annihilator $\sigma_n=\inf\{t\ge 0:
\rmK_n(t)=0\}$. We define further $K_{n,0}$ to be the number of
epochs $s_j<\sigma_n$ when none of the primary particles is marked
`head'. We define $U_n$ to be the number of particles that remain
at time $\sigma_n$; all these particles are secondary, so
$U_n=\Pi_n(\sigma_n)$. For two real-valued random variables $X$
and $Y$, we say that $X$ is stochastically smaller than $Y$ and
write $X\le_d Y$, if $P(X\le x)\ge P(Y\le x)$ for every
$x\in{\mathbb R}$.
\begin{lemma}\label{L1} The following stochastic order relations are true
   \begin{equation}\label{BD}
      K_n - K_{n,1}\ \le_d\ X_n\ \le_d\ K_n + K_{n,0}+X_{U_n},
   \end{equation}
   where in the right hand side $(K_n,K_{n,0},U_n)$ is assumed to be
   independent of $\{X_j:j\in\mn\}$.
\end{lemma}
\begin{proof}
We use the above coupling. By collision of at least two primary
particles both processes $\Pi_n$ and $\rm K_n$ decay, hence
$K_n-K_{n,1}\le X_n$. The number of jumps of $\Pi_n$ before and
including time $\sigma_n$ does not exceed the number of Poisson
epochs, which is $K_n+K_{n,0}$, and after time $\sigma_n$ the
coalescent evolves with $U_n$ particles, independently of the
history (given $U_n$). This yields the upper bound for $X_n$.
\end{proof}

\begin{lemma}\label{L2}
The sequence of distributions of $X_{U_n}$ $(n\in{\mathbb N})$ is
tight.
\end{lemma}
\begin{proof}
First of all, $1\leq X_{U_n}<U_n$ because the number of collisions
after $\sigma_n$ does not exceed the number of particles that
remain. To bound $U_n$, we focus on the behaviour of the secondary
particles.
At any time $s_j\leq \sigma_n$, each of the secondary particles is
removed with probability $1-\eta_j$ and at most one secondary
particle is born (possibly by collision of only primary
particles). We can only increase the number of secondary particles
if we assume that a secondary particle is always born at every
time $s_j\leq\sigma_n$.

For $j\in{\mathbb N}_0:={\mathbb N}\bigcup\{0\}$ denote by $Q_j$
the number of secondary particles at time $s_j$ ($s_0:=0$) and let
$\mathcal{F}$ be the $\sigma$-field generated by all the variables
$\eta_k$, $k\in{\mathbb N}$. Note that $Q_0=0$. With the above convention,
given ${\cal F}$, the sequence
$\{Q_j:j\in\mathbb{N}_0\}$ becomes a time-inhomogeneous Markov
chain on the set of non-negative integers with transition
probabilities
$$
\mmp\{Q_j=k\,|\,Q_{j-1}=m\} \ =\ {m\choose{k-1}} \eta_j^{k-1}
(1-\eta_j)^{m-k+1},\qquad k\in\{1,\ldots,m+1\}.
$$
Set
$$N_n:=\inf\{k\geq 1: n\eta_1\cdots\eta_k\leq 1\},\quad
n\in{\mathbb N}.$$ Almost surely on the event $\{N_n\geq 2\}$,
$$\mathbb{E}\left(Q_{N_n}|\mathcal{F}\right)=
1+\sum_{j=2}^{N_n}\eta_j\eta_{j+1}\cdots\eta_{N_n}=:1+\Theta_n, \
\ \ n\in{\mathbb N}.$$ Even though the subsequent argument is the
same as in the proof of \cite[Theorem 1.1]{IksTer}, we include it
for completeness. For $k=2,3,\ldots$, we have
$$\Theta_n1_{\{N_n=k\}}=\eta_1\ldots\eta_{N_n}\left({1\over\eta_1}+{1\over \eta_1\eta_2}+
\ldots+{1\over
\eta_1\ldots\eta_{N_n-1}}\right)
1_{\{N_n=k\}}\leq
$$
$$\leq {1\over
n}\left({1\over\eta_1}+{1\over \eta_1\eta_2}+ \ldots+{1\over
\eta_1\ldots\eta_{k-1}}\right)1_{\{N_n=k\}}\quad \text{a.s.},$$
which implies that $$\Theta_n1_{\{N_n\geq 2\}}\leq
{1\over n}\sum_{k=1}^{N_n-1}
\dfrac{1}{\eta_1\cdots\eta_k}
\quad
\text{a.s.}$$ By the key renewal theorem,
$$\lin {1\over n}\me \left( \sum_{k=0}^{N_n-1}
\dfrac{1}{\eta_1\cdots\eta_k}
\right)=$$$$=\lin
\int_0^{\log n}e^{-(\log
n-x)}\left(\sum_{k=0}^\infty\mmp\{-\log\eta_1-\ldots-\log\eta_k
\in {\rm d}x\}\right) =  \dfrac{1}{\me (-\log \eta)},$$ no matter
whether ${\tt m}_1:=\me (-\log\eta)$ is finite or infinite. Since
$\me Q_{N_n}1_{\{N_n=1\}}=o(1)$, we have proved that
$$\underset{n\to\infty}{\lim\sup}\,\me Q_{N_n}\leq 1/{\tt m}_1+1,$$
which implies that the sequence of distributions of $Q_{N_n}$
$(n\in{\mathbb N})$ is tight. Note that $N_n$ is the number of
steps that the multiplicative renewal process
$(\eta_1\cdots\eta_j,~j=0,1,\ldots)$, needs to drop below level
$1/n$. Now we wish to replace the fixed level $1/n$ by a random
level associated with the last primary particle disappearing at
time $\sigma_n$; this random level can be identified with the
minimal order statistic $R_{1,n}$ of the sample of size $n$ from
the uniform $[0,1]$ distribution (see the remark after Theorem
\ref{basic1} for this connection). From the subadditivity of the
first-passage time processes for random walks, and from the
convergence in distribution of $-\log R_{1,n}-\log n$ (which is a
standard fact from the extreme-value theory), it follows that
the number of steps $\rho_n$, say, the renewal process needs to
cross the interval $(\min(R_{1,n},1/n), \max(R_{1,n}, 1/n))$ remains
stochastically bounded as $n\to\infty$. On the event $A_n:=\{R_{1,n}\leq 1/n\}$,
as the process proceeds from $1/n$ to $R_{1,n}$, the number of
secondary particles may increase at most by $\rho_n$. On the
complementary event, $U_n-Q_{N_n}\leq \chi_n U_n$ a.s., where
$\chi_n\in (0,1)$ is a stochastic sequence (depending on $\rho_n$)
which is bounded away from one. The latter inequality follows from
the observation that in the coalescent process under consideration 
the number of particles involved in the generic collision is
proportional to the number of existing particles.
In conclusion, $|Q_{N_n}-U_n|$ is stochastically bounded on both
events, which completes the proof.
\end{proof}

\section{Number of collisions} \label{NrCol}

\par With these preliminaries, we are ready to derive the asymptotic
behaviour of the number of collisions $X_n$ from the previously
known results on $K_n$ \cite{GINR}. Introduce the logarithmic
moments
$${\tt m}_1:=\me(-\log\eta),\qquad {\tt m}_2:={\rm Var}(\log\eta)$$
for a random variable $1-\eta$ with distribution $\nu$ and recall
the standing assumptions of the paper given at the end of Section
\ref{intro}. It is worth mentioning explicitly that in parts
(1)-(3) of the next result ${\tt m}_1<\infty$.

\begin{thm}\label{basic1}
The following two assertions are equivalent:
\begin{enumerate}
\item[\rm (i)] There exist constants $\{a_n, b_n:n\in\mn\}$ with
$a_n>0$ and $b_n\in\mr$ such that, as $n\to\infty$,
$(X_n-b_n)/a_n$ converges weakly to some non-degenerate and proper
distribution.
\item[\rm (ii)] The distribution of $(-\log \eta)$
either belongs to the domain of attraction of a stable law, or the
function $x\mapsto\mmp(-\log \eta>x)$ slowly varies at $\infty$.
\end{enumerate}
Furthermore, there are five different regimes of convergence.
\begin{enumerate}
   \item[\rm(1)] If ${\tt m}_2<\infty$, then, with $b_n:={\tt m}_1^{-1}\log n$
      and $a_n:=({\tt m}_1^{-3}{\tt m}_2\log n)^{1/2}$, the limiting law
      of $(X_n-b_n)/a_n$ is standard normal.
   \item[\rm (2)] If ${\tt m}_2=\infty$ and
      $$
      \int_{(x,1)}(\log y)^2\,\nu({\rm d}y)
      \ \sim\ L(-\log x)
      \quad\text{as $x\to 0$},
      $$
      for some $L$ slowly varying at $\infty$, then, with
      $b_n:={\tt m}_1^{-1}\log n$ and $a_n:={\tt m}_1^{-3/2}c_{[\log n]}$,
      where $c_n$ is any sequence satisfying
      $\lim_{n\to\infty}nL(c_n)/c_n^2=1$, the limiting law of
      $(X_n-b_n)/a_n$ is standard normal.
   \item[\rm (3)] Assume that for some function $L$ slowly varying at
      $\infty$, the relation
      \begin{equation}\label{domain1}
         \mmp(\eta\leq x)\ \sim\ (-\log x)^{-\alpha}L(-\log x)
         \quad \text{as $x\to 0$},
      \end{equation}
      holds with $\alpha\in[1,2)$,  and that ${\tt m}_1<\infty$ if $\alpha=1$.
      Then, with $b_n:={\tt m}_1^{-1}\log n$ and
      $a_n:={\tt m}_1^{-(\alpha+1)/\alpha}c_{[\log n]}$, where $c_n$
      is any sequence satisfying $\lim_{n\to\infty}nL(c_n)/c_n^\alpha=1$,
      the limiting law of $(X_n-b_n)/a_n$ is $\alpha$-stable with
      characteristic function
      $$
      t\ \mapsto\ \exp\{-|t|^\alpha\Gamma(1-\alpha)(\cos(\pi\alpha/2)
      +i\sin(\pi\alpha/2)\, {\rm sgn}(t))\},\quad t\in\mr.
      $$
   \item[\rm (4)] Assume that ${\tt m}_1=\infty$ and that {\rm (\ref{domain1})}
      holds with $\alpha=1$. Let $c$ be any positive function satisfying
      $\lim_{x\to\infty}xL(c(x))/c(x)=1$ and set
      $\psi(x):=x\int_{\exp(-c(x))}^1 \mmp(\eta\leq y)/y\,dy$.
      Let $b(x)$ be any positive function satisfying
      $b(\psi(x))\sim\psi(b(x))\sim x$. Then, with $b_n:=b(\log n)$ and
      $a_n:=b(\log n)c(b(\log n))/\log n$, the limiting law of
      $(X_n-b_n)/a_n$ is $1$-stable with characteristic function
      $$
      t\ \mapsto\ \exp\{-|t|(\pi/2-i\log|t|\,{\rm sgn}(t))\},
      \quad t\in\mr.
      $$
   \item[\rm (5)] If {\rm (\ref{domain1})} holds with
      $\alpha\in [0,1)$, then, with $a_n:=\log^\alpha n/L(\log n)$,
      the limiting law of $X_n/a_n$ is scaled Mittag-Leffler
      $\theta_\alpha$  (exponential, if $\alpha=0$) with moments
      $$
      \int_{[0,\infty)} x^n\theta_\alpha({\rm d}x)
      \ =\ \dfrac{n!}{\Gamma^n(1-\alpha)\Gamma(1+n\alpha)},
      \quad n\in\mn.
      $$
\end{enumerate}
\end{thm}

\begin{proof}
Following the same line of argument as in the proof of Theorem 2.1
in \cite{GINR}, one can show that a necessary condition for the
convergence of $(X_n-b_n)/a_n$ is that the scaling satisfies
$a_n\to\infty$, which will be assumed in the sequel. To derive the
result from Theorem 2.3 in \cite{GINR} it is enough to  show that
$(X_n-b_n)/a_n$ weakly converges to a proper and non-degenerate
distribution if and only if $(K_n-b_n)/a_n$ weakly converges to
the same distribution. This equivalence is implied by (\ref{BD}).
Indeed, by Lemma \ref{L2}, $X_{U_n}/a_n$ converges in probability
to $0$. Under the assumption (\ref{nuat1}), $K_{n,0}$ converges in
distribution by \cite[Theorem 2.2]{GINR}, and $K_{n,1}$ converges
in distribution by \cite[Proposition 5.2]{GINR} (see also
\cite{small} for the case ${\tt m}_1<\infty$).
\end{proof}

\vskip0.1cm \noindent {\bf Remark.} We describe briefly a more
general construction with possibly infinite characteristic measure
$\nu$ satisfying $\int_{(0,1)}x\,\nu({\rm d}x)<\infty$. Let
$\tilde{\nu}$ be the image of $\nu$ under the transformation
$x\mapsto-\log x$, let $(\tilde{S}_t)$ be a subordinator with
L{\'e}vy measure $\tilde{\nu}$, and $S_t=\exp(-\tilde{S}_t)$. Then
$(S_t)$ is a non-increasing positive pure-jump process with
$S_0=1$. Primary particles are labelled by independent
uniform$[0,1]$ points $R_1,\ldots,R_n$, also independent of
$(S_t)$. The annihilator process ${\rm K}_n=({\rm K}_n(t), t\geq
0)$ is defined by setting ${\rm K}_n(0)=n$ and by decrementing the
process at time $t$ by $k$ if the interval $[S_{t}, S_{t-}]$
contains exactly $k$ of the uniform points. Thus ${\rm K}_n$ may
have a jump at $t$ only if $S_{t-}>S_t$. The coalescent $\Pi_n$ is
defined in a similar way, but except $R_1,\ldots,R_n$ representing
primary particles, new points representing secondary particles
appear. Set $\Pi_n(0)=n$. Conditionally given that $[S_{t},
S_{t-}]$ contains $k\geq 1$ points, primary and secondary, let
$\Pi_n$ decrement by $k-1$ and let one point with uniform
distribution on $[0,S_t]$ appear to represent a secondary
particle. The inequality $X_n\geq K_n-K_{n,1}$ is still true, but
another bound in (\ref{BD}) makes no sense for infinite $\nu$,
since a variable  analogous to $K_{n,0}$ is identically $\infty$.
It seems plausible that in this more general setting $X_n$ and
$K_n-K_{n,1}$ still have the same limit laws and require the same
scaling. Asymptotics of $K_n, K_{n,1}$ were studied for infinite
measures $\nu$ satisfying a condition of regular variation
\cite{GnePitYor1, GPYI}.

\begin{example} \label{beta}
\rm Application of case (1) of Theorem \ref{basic1} to coalescents
governed by beta measures $\nu({\rm d}x)=c x^{a-3}(1-x)^{b-1}{\rm
d}x$ with parameters $a>2$ and $b,c>0$ shows that the scaled and
centered number of collisions
$$(X_n-{\tt m}_1^{-1}\log n)/ ({\tt m}_2{\tt
m}_1^{-3}\log n)^{1/2}$$ converges in distribution to the standard
normal law, where ${\tt m}_1=\Psi(a-2+b)-\Psi(b)$ and ${\tt
m}_2=\Psi^\prime(b)-\Psi^\prime(a-2+b)$, and $\Psi$ is the
logarithmic derivative of the gamma function.

\par For the special case $\nu({\rm d}x)={\rm d}x$ (uniform distribution)
we have $g_{n,m}/g_n=1/(n-1)$. Thus, $\Pi_n$ jumps from $n$ to
every $m\in\{1,\ldots,n-1\}$ with the same probability. In this
case Theorem \ref{basic1} implies that $(X_n-\log n)/\sqrt{\log
n}$ is asymptotically normal, which agrees with a classical result
known in the context of records or cyclic partition of
permutations. We refer to the survey of \cite{Bai} for some
further information along with references.

\par More generally, when $\nu$ is a beta$(1,b)$ distribution, a
straightforward calculation shows that $\Pi_n$ jumps from $n\geq
2$ to any state $m\in\{1,\ldots,n-1\}$ with probability
$w_m/(w_1+\cdots+w_{n-1})$, where $w_k:=\Gamma(k+b-1)/\Gamma(k)$.
By \cite[Proposition 2]{Ross}, $X_n$ is representable as the sum
of $n-1$ independent indicators $I_k$ with
$$
\mmp(I_k=1)\ =\ \dfrac{w_k}{w_1+\cdots+w_k},\qquad
k\in\{1,\ldots,n-1\}.
$$
It can be shown that $\mmp(I_k=1)\sim b/k$ as $k\to\infty$. From
this, as $n\to\infty$, the expectation and the variance of $X_n$
are both asymptotically equal to $b\log n$. By the central limit
theorem, the distribution of $(X_n-b\log n)/\sqrt{b\log n}$
converges to the standard normal law, which again agrees with
Theorem \ref{basic1}, case (1). The latter result is classics
which dates back to \cite{Gon}.
\end{example}

\par The merging evolution of $n$ particles in the coalescent is
often represented by a rooted tree with $n$ leaves, and with the
internal nodes corresponding to collisions \cite{Bertoin, CSP}.
One interpretation of the similarity of the asymptotic behaviour
of $X_n$ and $K_n$ is that almost all internal nodes of the tree
are directly connected to one of the leaves.

\par Some other properties of the coalescent follow easily in the same way.
Let $Z_n$ be the time of the first collision of a tagged particle
(say, labelled 1). In the tree representation, $Z_n$ is the length
of an external branch, which connects a leaf to its closest
internal node. The asymptotic behaviour of $Z_n$ for the Kingman
coalescent and the Bolthausen-Sznitman coalescent ($\nu({\rm
d}x)=x^{-2}{\rm d}x$) was studied in \cite{Amke} and
\cite{FreMoe}, respectively. In our setting we have the following.

\begin{assertion} As $n\to\infty$, $Z_n$ converges weakly to an
exponential variable with parameter $p:={\mathbb
E}(1-\eta)=\int_{(0,1)} x\nu({\rm d}x)$. Moreover, the number of
collisions in the coalescent before and including time $Z_n$
converges weakly to a geometric variable with parameter $p$.
\end{assertion}
\begin{proof}
In the process with infinitely many particles a tagged particle is
marked `head' at rate $p$, and when this occurs the particle
collides with infinitely many other particles, because the
assumption that $\nu$ is a probability measure implies that the
coalescent does not come down from infinity \cite{Pit} (which is
also obvious from the construction). This proves the first claim.
The second claim follows similarly, by noting that the time of the
first collision of a given primary particle is $s_1$ with
probability $p$.
\end{proof}

\section{Absorption times} \label{AbsTimes}

The coupling allows us to derive the asymptotics of $\tau_n$.
\begin{lemma}\label{L3}
The following stochastic order relations hold
$$\sigma_n\ \le_d\ \tau_n\ \le_d\ \sigma_n+\tau'_{U_n},$$
where, on the right-hand side, $\tau_j'\stackrel{d}{=}\tau_j$,
$j\in{\mathbb N}$, and $\sigma_n,(\tau_j')$ and $U_n$ are assumed
to be independent.
\end{lemma}
\begin{proof}
These relations are obvious from the coupling. When all the
primary particles disappear at time $\sigma_n$, there are $U_n$
secondary particles left.
\end{proof}

From the results obtained in \cite{GINR} it follows that
$\sigma_n$ has the same limiting distribution as the first passage
time through the level $\log n$ for a compound Poisson process
$\{R_t: t\geq 0\}$ with the generic jump $(-\log \eta)$ and
intensity one. By Lemma \ref{L2} and Lemma \ref{L3}, the same
asymptotics is valid for $\tau_n$. This argument leads to the
following result.

\begin{thm}\label{main1}
Assume that condition {\rm (ii)} of Theorem \ref{basic1} holds.
Then, with $X_n$ replaced by $\tau_n$, condition {\rm (i)} and
parts {\rm(2) - (5)} of that theorem hold without changes, while
part {\rm (1)} holds with $a_n=({\tt m}_1^{-3}({\tt m}_2+{\tt
m}_1^2)\log n)^{1/2}$.

\end{thm}
\begin{proof}
We focus on the asymptotics of the first passage time process
$\{U_x: x\geq 0\}$ defined as
$$U_x:=\inf\{t\ge 0:R_t\geq x\},\qquad x\ge 0.$$
Let $N_x$ be the number of renewals on $[0,x]$ of the renewal
process with $N_0=0$ and steps distributed as $(-\log \eta)$.
Obviously, $U_x=s_{N_x}$, where as above $s_1,s_2,\ldots$ are the
arrival times of an independent unit Poisson process. Assume that
condition (ii) of Theorem \ref{basic1} holds, and let $a(x)$ and
$b(x)$ be the same as in \cite[Proposition 8.1]{GINR} (setting
$a_n:=a(\log n)$ and $b_n:=b(\log n)$ gives the normalizing
constants for the cases (2) - (5) of Theorem \ref{basic1}). We
will exploit the following decomposition
\begin{equation}\label{kr}
\dfrac{s_{N_x}-b(x)}{a(x)}=\dfrac{s_{N_x}-N_x}{\sqrt{N_x}}\dfrac{\sqrt{N_x}}{a(x)}+\dfrac{N_x-b(x)}{a(x)}.
\end{equation}
In view of the independence and since $N_x\to\infty$ almost surely
as $x\to\infty$, we have by the central limit theorem that, as
$x\to\infty$, $(s_{N_x}-N_x)/\sqrt{N_x}$ converges in distribution
to a random variable $\xi_1$ with the standard normal
distribution. By \cite[Proposition 8.1]{GINR}, the right-most term
in (\ref{kr}) converges in distribution to a random variable
$\xi_2$ with a proper and non-degenerate law. Note that $\xi_1$
and $\xi_2$ are independent.

\par Case (1). By \cite[Proposition 8.1]{GINR}, $b(x)={\tt m}_1^{-1}x$,
$a(x)=({\tt m}_2{\tt m}_1^{-3}x)^{1/2}$, and $\xi_2$ has the
standard normal law. Since, by the strong law of large numbers, as
$x\to\infty$, $N_x/x\to {\tt m}_1^{-1}$ almost surely, we conclude
that the left-hand side in (\ref{kr}) converges in distribution to
a random variable ${\tt m}_1 {\tt m}_2^{-1/2}\xi_1+\xi_2$ which
has the normal law with zero mean and variance ${\tt m}_1^2 {\tt
m}_2^{-1}+1$. Therefore, $(U_x-{\tt m}_1^{-1}x)/({\tt
m}_1^{-3}({\tt m}_2+{\tt m}_1^2)x)^{1/2}$ weakly converges to the
standard normal law, as desired.
\par Cases (2) - (5). By \cite[Proposition 8.1]{GINR}, as
$x\to\infty$, $\sqrt{N_x}/a(x)$ converges to $0$ in probability.
Therefore, the left-hand side in (\ref{kr}) weakly converges to
the same law as the right-most term in (\ref{kr}). We again refer
to \cite[Proposition 8.1]{GINR} for the identification of the
normalizing constants and the relevant limiting laws.
\end{proof}

\begin{example}
\rm Theorem 3.2 is obviously applicable to the class of
beta-coalescents with parameter $a>2$ and $b,c>0$ considered in
Example \ref{beta}. Consequently, for such beta-coalescents, the
scaled and centered absorption time $(\tau_n-{\tt m}_1^{-1}\log
n)/({\tt m}_1^{-3}({\tt m}_2+{\tt m}_1^2)\log n)^{1/2}$ converges
in distribution to the standard normal law, with ${\tt m}_1$ and
${\tt m}_2$ as given in Example \ref{beta}.
\end{example}

\vskip0.1cm \noindent {\bf Acknowledgement} A. Iksanov gratefully
acknowledges the support by the German Scientific Foundation
(project no. 436UKR 113/93/0-1).


\begin{thebibliography}{99}
\footnotesize

\bibitem{Bai}{\sc Bai, Z.D., Hwang, H.K., and Liang, W.Q.} (1998).
Normal approximations of the number of records in geometrically
distributed random variables. {\em Random Struct. Algorithms} {\bf
13}, 319--334.

\bibitem{BBS} {\sc Berestycki, J., Berestycki, N. and  Schweinsberg, J.} (2008)
Small time properties of beta coalescents, 
 {\em Ann. Inst. H. Poincar{\'e} Probab. Statist.}   {\bf 44}, 214-238. 


\bibitem{Bertoin}{\sc Bertoin, J.} (2006). {\it Random fragmentation
and coagulation processes}, Cambridge Univ. Press.

\bibitem{Amke}{\sc Caliebe, A., Neininger, R., Krawczak, M., and R\"osler,
U.} (2007). On the length distribution of external branches in
coalescence trees: genetic diversity within species. {\em Theor.
Popul. Biol.} {\bf 72}, 245–-252.

\bibitem{Delmas}{\sc Delmas, J.-F., Dhersin, J.-S. and Siri-Jegousse, A.}
(2008). Asymptotic results on the length of coalescent trees. {\em
Ann. Appl. Probab.} {\bf 18,} 997–-1025.

\bibitem{DrmIksMoeRoe}{\sc Drmota, M., Iksanov, A., M\"ohle, M., and R\"osler, U.}
(2007). Asymptotic results concerning the total branch length of
the Bolthausen-Sznitman coalescent. {\em Stoch. Process. Appl.}
{\bf 117}, 1404--1421.

\bibitem{DrmIksMoeRoe2}{\sc Drmota, M., Iksanov, A., M\"ohle, M.,
and R\"osler, U.} (2008). A limiting distribution for the number
of cuts needed to isolate the root of a random recursive tree.
{\em Random Struct. Algorthms}, to appear.

\bibitem{FreMoe} {\sc Freund, F. and M\"ohle, M.} (2007).
On the time back to the most recent common ancestor and the
external branch length of the Bolthausen-Sznitman coalescent. {\em
Markov Process. Relat. Fields}, submitted.

\bibitem{Gne}{\sc Gnedin, A.V.} (2004). The Bernoulli sieve.
{\em Bernoulli} {\bf 10}, 79--96.

\bibitem{GDP} {\sc Gnedin, A., Dong, R., and  Pitman, J. } (2007).
Excheangable partitions derived from Markovian coalescents. {\it
Ann. Appl. Probab.} {\bf  17}, 1172--1201.

\bibitem{small}{\sc Gnedin, A., Iksanov, A., and Roesler, U.} (2008).
Small parts in the Bernoulli sieve, available from arXiv:math.

\bibitem{GINR}{\sc Gnedin, A., Iksanov, A., Negadajlov, P., and Roesler, U.} (2008).
The Bernoulli sieve revisited, available from arXiv:math.

\bibitem{GnePitYor1}{\sc Gnedin, A., Pitman, J., and Yor, M.} (2006).
Asymptotic laws for regenerative compositions: gamma subordinators
and the like. {\em Probab. Theory Relat. Fields} {\bf 135},
576--602.

\bibitem{GPYI}{\sc Gnedin, A., Pitman, J., and Yor, M.} (2006).
Asymptotic laws for compositions derived from transformed
subordinators. {\em Ann. Probab.} {\bf 34}, 468--492.

\bibitem{RCS} {\sc Gnedin, A. and Pitman, J.} (2005).
Regenerative composition structures. {\em  Ann. Probab. } {\bf
33}, 445--479.


\bibitem{GneYak}{\sc Gnedin, A. and Yakubovich, Y.} (2007).
On the number of collisions in $\Lambda$-coalescents. {\em
Electron. J. Probab.} {\bf 12}, 1547--1567.

\bibitem{Gon} {\sc Goncharov, V.L.} (1962). On the field of combinatorial analysis.
{\em Amer. Math. Soc. Translations}, series 2, {\bf 19}, 1--46.

\bibitem{IksMarMoe}{\sc Iksanov, A., Marynych, A., and M\"ohle, M.} (2007).
On the number of collisions in beta(2,b)-coalescents, submitted.

\bibitem{IksMoe2}{\sc Iksanov, A. and M\"ohle, M.} (2008).
On the number of jumps of random walks with a barrier. {\em Adv.
Appl. Probab.} {\bf 40}, 206--228.

\bibitem{IksTer}{\sc Iksanov, A. and Terletsky, Yu.} (2008). On
asymptotic behavior of certain recursions with random indices of
linear growth, submitted.

\bibitem{Pan}{\sc Panholzer, A.} (2004). Desctruction of recursive
trees, in: {\em Mathematics and Computer Science} III,
Birkh\"auser, Basel, pp. 267--280.

\bibitem{Pit}{\sc Pitman, J.} (1999). Coalescents with multiple
collisions. {\em Ann. Probab.} {\bf 27}, 1870--1902.

\bibitem{CSP}{\sc Pitman, J.} (2006). {\it Combinatorial stochastic processes}.
Lecture Notes in Mathematics { \bf 1875}, Springer.

\bibitem{Ross}{\sc Ross, S.M.} (1982). A simple heuristic approach to
simplex efficiency. {\em Europ. J. Operat. Research} {\bf 9},
344--346.

\bibitem{Sag}{\sc Sagitov, S.} (1999). The general coalescent with asynchronous
mergers of ancestral lines. {\em J. Appl. Probab.} {\bf 36}, 1116--1125.

\end{thebibliography}
\end{document}